\documentclass[12pt,reqno]{amsart}
\usepackage{fullpage}
\usepackage[top=1.5in, bottom=2in, left=2.5in, right=2.5in]{geometry}
\usepackage{graphicx,amsmath,amssymb,amsfonts,amsthm,enumitem,bm,xcolor}
\usepackage{fullpage}
\usepackage{cleveref}

\usepackage[normalem]{ulem}

\usepackage{url}

\usepackage{color}

\newcommand{\z}{\zeta}

\renewcommand{\(}{\left\(}
\renewcommand{\)}{\right\)}

\numberwithin{equation}{section}
\theoremstyle{plain}
\newtheorem{theorem}{Theorem}[section]
\newtheorem{lemma}[theorem]{Lemma}

\newtheorem{conjecture}[theorem]{Conjecture}
\newtheorem*{conjecture*}{Conjecture}

\theoremstyle{definition}

\numberwithin{equation}{section}

\renewcommand{\pmod}[1]{\ \left( \mathrm{mod} \, #1 \right)}

\setlist[enumerate]{leftmargin=*,label=\rm{(\arabic*)}}

\makeatletter
\@namedef{subjclassname@2020}{%
	\textup{2020} Mathematics Subject Classification}
\makeatother

\allowdisplaybreaks

\title[]{Partition statistics and universal mixed mock modular forms}
\author{Koustav Banerjee}
\author{Kathrin Bringmann}
\address{University of Cologne, Department of Mathematics and Computer Science, Weyertal 86-90, 50931 Cologne, Germany}
\email{kbanerj1@uni-koeln.de}
\email{kbringma@math.uni-koeln.de}

\makeatletter
\@namedef{subjclassname@2020}{%
	\textup{2020} Mathematics Subject Classification}
\makeatother


\begin{document}

\begin{abstract}
 In this note, we construct combinatorial statistics in terms of partitions, overpartitions, and two-color partitions and relate the corresponding generating functions to the universal mock theta functions. Our framework primarily builds on two $q$-series identities of Fine. 
\end{abstract}

\maketitle

\section{Introduction and Statements of results} 

 In his last letter to Hardy, Ramanujan \cite[pp. 354--355]{BB} introduced seventeen $q$-series, which he called  mock theta functions. According to Ramanujan, a {\it mock theta function} is a function $F$ whose asymptotic behavior near roots of unity resembles that of weakly holomorphic modular forms, however, there does not exist a single weakly holomorphic modular form whose difference from $F$ remains bounded as $q$ approaches any root of unity radially. In his seminal thesis, Zwegers \cite{Zw} explained the modularity of Ramanujan's mock theta functions by ``completing" them to so-called harmonic Maass forms.\footnote{A {\it harmonic Maass form} is a function that satisfies modularity, is annihilated by the weight $k$ hyperbolic Laplacian, and satisfies a growth condition at the cusps. Its holomorphic part is called a {\it mock modular form}.} 
 
 Next recall that Dyson \cite{Dyson} introduced the {\it rank}\,\footnote{The {\it rank} of a partition is defined as its largest part minus the number of its parts.} statistic to give a combinatorial explanation of partition congruences of Ramanujan (see  \cite{Ramanujan}). Note that the generating function of the number of partitions with even rank minus the number of partitions with odd rank is Ramanujan's third order mock theta function \[f(q):=\sum_{n\ge 0}\frac{q^{n^2}}{(-q)_n^2}.\]
 Here, and throughout, $|q|<1$ and for $a\in\mathbb{C}$ and $n\in\mathbb{N}_0\cup \{\infty\}$, we set $(a)_n\hspace{-0.1 cm}=\hspace{-0.1 cm}(a;q)_n\hspace{-0.1 cm}:=\prod_{j=0}^{n-1}(1-aq^j)$. 
 
 Moreover,  MacMahon \cite{M} studied partitions without sequences which are related to mixed mock modular forms.\footnote{Roughly speaking, a {\it mixed mock modular form} is a linear combination of products of modular forms and mock modular forms.} More precisely, using work of Fine \cite{Fine}, Andrews \cite{Andrews} showed that the generating function of partitions without sequences is related to Ramanujan's third order mock theta function $\chi(q)$. Note that Ramanujan's mock theta functions are specializations of the {\it universal mock theta functions} (see \cite[(3.2) and above (5.1)]{GM})
\begin{equation}\label{unimockdef}
g_3(\zeta;q):=\sum_{n\ge 0}\frac{q^{n(n+1)}}{\left(\zeta,\frac q\zeta\right)_{n+1}},\qquad g_2(\zeta;q):=\sum_{n\ge 0}\frac{(-q)_n }{\left(\zeta,\frac q\zeta\right)_{n+1}}q^{\frac{n(n+1)}{2}}.
\end{equation}
 Here $\z\in \mathbb{C}\setminus\{0\}$ and for $r\in \mathbb{N}$ and $a_1,\dots,a_r\in \mathbb{C}$, we let $(a_1,\dots,a_r)_n:=(a_1,\dots,a_r;q)_n:=(a_1)_n\dots (a_r)_n$.

In \cite[Theorem 1.3]{BLM}, Lovejoy, Mahlburg, and the second author gave a combinatorial interpretation of (for $d\ge3$ and $d>2r$)
\[
\frac{\left(-q^r,-q^{d-r};q^d\right)_\infty}{\left(q^{2d};q^{2d}\right)_\infty}g_2\left(-q^r;q^d\right).
\]
Moreover, \cite[(1.2)]{BLM} yields a combinatorial statistic that is related to (for $d\ge3$ and $d>2r$)
\[
\left(-q^r,-q^{d-r};q^d\right)_\infty g_3\left(-q^r;q^d\right).
\]
Note that $d=2r$ was excluded in \cite{BLM}. In this case, the residue classes $r$ and $-r\pmod{d}$ coincide. Thus $(-q^r,-q^{d-r};q^d)_\infty=(-q^r;q^{2r})_\infty^2$ and each part $\equiv r\pmod{2r}$ comes in two copies. This naturally leads to overpartitions and two-color partitions. This is a key motivation for our combinatorial framework. In this note, we consider restricted partitions, overpartitions, and colored partitions and relate the corresponding generating functions to the universal mock theta functions. We begin with different $q$-hypergeometric  series representations of $g_3$ and $g_2$. To the authors' knowledge these representations are not known. 

\begin{theorem}\label{sec:ummmflem1}
Let $\z\neq q^{m}$ for any $m\in \mathbb{Z}$.
\begin{enumerate}[label=\rm(\arabic*)]
		\item We have
		\begin{equation*}%\label{sec:ummmflem1eqn1}
		g_3(\zeta;q)=\frac{(q)_{\infty}}{\left(\zeta,\frac q\zeta\right)_{\infty}}\sum_{n\ge 0}\frac{\left(\zeta,\frac q\zeta\right)_n}{(q)_n}q^n.
		\end{equation*}
		\item We have
		\begin{equation*}%\label{sec:ummmflem1eqn2}
		g_2(\zeta;q)=\frac{\left(q^2;q^2\right)_{\infty}}{\left(\zeta,\frac q\zeta\right)_{\infty}}\sum_{n\ge 0}\frac{\left(\zeta,\frac q\zeta\right)_n}{\left(q^2;q^2\right)_n}q^n.
		\end{equation*}
	\end{enumerate}
\end{theorem}

With a suitable specialization of the variables $\z$ and $q$ and after multiplying by appropriate infinite products, the identities in \Cref{sec:ummmflem1} yield the generating functions for classes of restricted partitions. The summation index records either the largest part in a certain residue class or the number of parts in that class, while the remaining factors impose restrictions on the parts in other residue classes.

We now describe these partition statistics explicitly. Here and throughout, we assume $d,r\in \mathbb{N}$ with $d\ge 2r$. First, recall that an {\it overpartition} of $n$ is a partition of $n$ in which the first occurrence of a part may be overlined. Let $\overline{p}(n)$ denote the number of overpartitions of $n$ and set $\overline{p}(0):=1$. For instance, $\overline{p}(2)=4$ enumerated by $2, \overline{2}, 1+1, \overline{1}+1$. We now consider overpartitions of $n$ with the following restrictions:
\begin{enumerate}
	\item parts are divisible by $r$,
	\item parts that are divisible by $2r$ are non-overlined,
	\item the largest part is divisible by $2r$,
	\item each part $\equiv  r\pmod{2r}$ can appear at most twice and if it appears twice, then one of the two occurrences must be overlined.
\end{enumerate}
Let $\mathrm{op}^{[1]}_{r}(n)$ denote the number of such overpartitions of $n$ and set  $\mathrm{op}^{[1]}_{r}(0):=1$. For example, $\mathrm{op}^{[1]}_2(12)=5$ enumerated by $12$, $8+4$, $8+\overline{2}+2$, $4+4+4$, $4+4+\overline{2}+2$. Let $\operatorname{OP}^{[1]}_{r}(q):=\sum_{n\ge 0}{\rm op}^{[1]}_r(n)q^n$ denote the corresponding generating function.

\begin{theorem}\label{mmmfocase1lem2}
We have
	\[
	\operatorname{OP}^{[1]}_{r}(q)=\frac{\left(-q^r;q^{2r}\right)^2_{\infty}}{\left(q^{2r};q^{2r}\right)_{\infty}} g_3\left(-q^r;q^{2r}\right).
	\]
\end{theorem}

Next, for $d>2r$, we consider partitions of $n$ with the following restrictions:
\begin{enumerate}
	\item parts are $\equiv 0,\pm r\ (\text{mod}\ d)$,
    \item each part $\equiv \pm r\ (\text{mod}\ d)$ is greater than every part divisible by $d$.
\end{enumerate}
Let $p^{[2]}_{d,r}(n)$ denote the number of such partitions of $n$ and set $p^{[2]}_{d,r}(0):=1$. For example, $p^{[2]}_{4,1}(5)=3$ enumerated by $5$, $3+1+1$, $1+1+1+1+1$. Let $P^{[2]}_{d,r}(q):=\sum_{n\ge 0}p^{[2]}_{d,r}(n)q^n$ denote the corresponding generating function.

\begin{theorem}\label{mmmfocase2lem1}
We have, for $d>2r$,
	\[
	P^{[2]}_{d,r}(q)=\frac{1}{\left(q^d;q^d\right)_{\infty}} g_3\left(q^r;q^d\right).
	\]
\end{theorem}

 Next, recall that a {\it two-color partition} of $n$ is a partition of $n$ in which each part appears in one of two fixed colors, say blue or green. Parts of the same size but with different colors are regarded as distinct and if a part occurs both in blue and in green, then we list the green parts first. Let $p_2(n)$ denote the number of two-color partitions of $n$ and set $p_2(0):=1$. For example, $p_2(2)=5$ enumerated by $2_{g}$, $2_{b}$, $1_{g}+1_{g}$, $1_{g}+1_{b}$, $1_{b}+1_{b}$. We consider two-color partitions of $n$ with the following restrictions:
\begin{enumerate}
	\item parts are divisible by $r$,
	\item parts divisible by $2r$ can appear only in blue,
	\item each part $\equiv r\ (\text{mod}\ 2r)$ is greater than parts divisible by $2r$.
\end{enumerate} 
Let $p^{[2]}_{r}(n)$ denote the number of such two-color partitions of $n$ and set $p^{[2]}_{r}(0):=1$. For example, $p^{[2]}_2(8)=11$ enumerated by $8_b$, $6_g+2_g$, $6_b+2_g$, $6_g+2_b$, $6_b+2_b$, $4_b+4_b$, $2_g+2_g+2_g+2_g$, $2_g+2_g+2_g+2_b$, $2_g+2_g+2_b+2_b$, $2_g+2_b+2_b+2_b$, $2_b+2_b+2_b+2_b$. Let $P^{[2]}_{r}(q):=\sum_{n\ge 0}p^{[2]}_{r}(n)q^n$ denote the corresponding generating function.
\begin{theorem}\label{mmmfocase2lem2}
We have
	\[
	P^{[2]}_{r}(q)=\frac{1}{\left(q^{2r};q^{2r}\right)_{\infty}} g_3\left(q^r;q^{2r}\right).
	\]
\end{theorem}

Now, we consider overpartitions of $n$ with the following restrictions:
\begin{enumerate}
\item parts are divisible by $r$ but not by $4r$,
\item	parts $\equiv 2r\pmod{4r}$ are non-overlined,
\item	each part $\equiv r\pmod{2r}$ can appear at most twice, and if it appears twice, it appears once overlined and once non-overlined,
\item if there are exactly $m$ parts $\equiv 2r\pmod{4r}$ (counted with multiplicity), then every part $\equiv r\pmod{2r}$ is less than $2rm$. In particular, if there is no part $\equiv 2r \pmod{4r}$, then there is no part $\equiv r \pmod{2r}$.
\end{enumerate}
Let $\mathrm{op}^{[2]}_{r}(n)$ denote the number of such overpartitions of $n$ and set  $\mathrm{op}^{[2]}_{r}(0):=1$. For example, ${\rm op}^{[2]}_1(6)=3$ enumerated by $6$, $2+2+2$, $2+2+\overline{1}+1$. Let $\operatorname{OP}^{[2]}_{r}(q):=\sum_{n\ge 0}\mathrm{op}^{[2]}_{r}(n)q^n$ denote the corresponding generating function.
\begin{theorem}\label{mmmfecase1lem2}
	We have
	\[
	\operatorname{OP}^{[2]}_{r}(q)=\frac{\left(-q^r;q^{2r}\right)^2_{\infty}}{\left(q^{4r};q^{4r}\right)_{\infty}} g_2\left(-q^r;q^{2r}\right).
	\]
\end{theorem}

Next, for $d>2r$, we consider partitions of $n$ with the following restrictions: 
\begin{enumerate}[label=\rm(\arabic*)]
	\item parts are $\equiv d\!\pmod{2d}$ or $\equiv \pm r\!\pmod d$,
	\item if there are exactly $m$ parts $\equiv d\!\pmod{2d}$ (counted with multiplicity), then each part $\equiv \pm r\!\pmod d$ is larger than $dm$.
\end{enumerate}
Let $p^{[4]}_{d,r}(n)$ denote the number of such partitions of $n$ and set $p^{[4]}_{d,r}(0):=1$. For example, $p^{[4]}_{4,1}(5)=3$ enumerated by $5$, $3+1+1$, $1+1+1+1+1$. Let $P^{[4]}_{d,r}(q):=\sum_{n\ge 0}p^{[4]}_{d,r}(n)q^n$ denote the corresponding generating function.
\begin{theorem}\label{mmmfecase2lem1}
We have, for $d>2r$,
	\[
	P^{[4]}_{d,r}(q)=\frac{1}{\left(q^{2d};q^{2d}\right)_{\infty}} g_2\left(q^r;q^d\right).
	\]
\end{theorem}

Finally, for $d=2r$, we consider two-color partitions of $n$ with the following restrictions:
\begin{enumerate}[label=\rm(\arabic*)]
	\item parts are divisible by $r$ but not by $4r$,
	\item parts $\equiv 2r\!\pmod{4r}$ appear only in blue,
	\item if there are exactly $m$ parts $\equiv 2r\!\pmod{4r}$ (counted with multiplicity), then every part $\equiv r\!\pmod{2r}$ is larger than $2rm$.
\end{enumerate}
Let $p^{[3]}_r(n)$ denote the number of such two-color partitions of $n$ and set $p^{[3]}_r(0):=1$. For example, $p^{[3]}_2(8)=10$ enumerated by $6_g+2_g$, $6_b+2_g$, $6_g+2_b$, $6_b+2_b$, $4_b+4_b$, $2_g+2_g+2_g+2_g$, $2_g+2_g+2_g+2_b$, $2_g+2_g+2_b+2_b$, $2_g+2_b+2_b+2_b$, $2_b+2_b+2_b+2_b$. Let $P^{[3]}_r(q):=\sum_{n\ge 0}p^{[3]}_r(n)q^n$ denote the corresponding generating function.
\begin{theorem}\label{mmmfecase2lem2}
We have
	\[
	P^{[3]}_{r}(q)=\frac{1}{\left(q^{4r};q^{4r}\right)_{\infty}} g_2\left(q^r;q^{2r}\right).
	\]
\end{theorem}

This note is organized as follows. In \Cref{s1}, we recall $q$-series transformations which are required to prove \Cref{sec:ummmflem1}. \Cref{s2} presents the proofs of Theorems \ref{sec:ummmflem1}--\ref{mmmfecase2lem2}. \Cref{prob} is devoted to open problems for future research.

\section*{Acknowledgements}
The authors have received funding from the European Research Council (ERC)
under the European Union’s Horizon 2020 research and innovation programme (grant agreement
No. 101001179).

\section{Preliminaries}\label{s1}
In this section we recall $q$-series identities\footnote{Throughout, we use $q$-series identities whenever both sides are well-defined.} which are used in the proofs of our results. We start with Fine's identities. To state these, following \cite[equation (1.1)]{Fine}, for $a\in \mathbb{C}$, $b\neq q^{-m}$ for every $m\in \mathbb{N}$, and $|w|<1$, we define
\begin{equation*}%\label{sec0alem1eqn1a}
F(a,b;w;q):=\sum_{n\ge 0}\frac{(aq)_n}{(bq)_n}w^n.
\end{equation*}
By \cite[equation (11.3)]{Fine}, we have the following identity. 
\begin{lemma}\label{sec0blemma2eqn2}
	We have
	\[
	\frac{\left(bq, wq\right)_{\infty}}{\left(q\right)_{\infty}}F(0,b;w;q)=\frac{1}{1-w}\sum_{n\ge0}\frac{\left(b,w\right)_n}{\left(q\right)_n}q^n.
	\]
\end{lemma}

Next, from \cite[(12.3)]{Fine}, we obtain the following lemma.
\begin{lemma}\label{sec0blemma2eqn3}
	We have
	\[
	(1-w)F(0,b;w;q)=\sum_{n\ge0}\frac{(bw)^n }{\left(bq, wq\right)_n}q^{n^2}.
	\] 
\end{lemma}

Next we require the following transformation from \cite[(III.10)]{GR}.
\begin{lemma}\label{sec:ummmflem1eqn5}
We have
\begin{equation*}
\sum_{n\ge 0}\frac{\left(\frac{aq}{bc},d,e\right)_n}{\left(q,\frac{aq}{b},\frac{aq}{c}\right)_n}\left(\frac{aq}{de}\right)^n%\\
=\frac{\left(\frac{aq}{d},\frac{aq}{e},\frac{aq}{bc}\right)_{\infty}}{\left(\frac{aq}{b},\frac{aq}{c},\frac{aq}{de}\right)_{\infty}}\sum_{n\ge 0}\frac{\left(\frac{aq}{de},b,c\right)_n}{\left(q,\frac{aq}{d},\frac{aq}{e}\right)_n}\left(\frac{aq}{bc}\right)^n.
\end{equation*}
\end{lemma}

Moreover, we require the following identity due to Euler (see \cite[(2.2.5)]{AndrewsBook}).
\begin{lemma}\label{Euler}
 For $|w|,|q|<1$, we have
\[
\sum_{n\ge0}\frac{w^n}{(q)_n}=\frac{1}{(w)_\infty}.
\]
\end{lemma}

\section{Proof of Theorems \ref{sec:ummmflem1}--\ref{mmmfecase2lem2}}\label{s2}

We start with the proof of \Cref{sec:ummmflem1}.

\begin{proof}[Proof of \Cref{sec:ummmflem1}]
	Combining Lemmas \ref{sec0blemma2eqn2} and \ref{sec0blemma2eqn3}, we obtain
	\begin{equation*}
	\sum_{n\ge0}\frac{(bw)^n }{\left(bq, wq\right)_n}q^{n^2}=\frac{(q)_{\infty}}{\left(bq,wq\right)_{\infty}}\sum_{n\ge0}\frac{\left(b, w\right)_n}{\left(q\right)_n}q^n.
	\end{equation*}
	Letting $(b,w)=(\zeta,\frac q\zeta)$, it follows that
	\begin{equation*}%\label{sec:ummmflem1eqn3}
	\sum_{n\ge 0}\frac{q^{n(n+1)}}{\left(\zeta q,\frac {q^2}\zeta\right)_{n}}=\frac{(q)_{\infty}}{\left(\zeta q,\frac{q^2} \zeta\right)_{\infty}}\sum_{n\ge0}\frac{\left(\zeta,\frac q\zeta\right)_n}{\left(q\right)_n}q^n.
	\end{equation*}
 This and the definition of $g_3$ from \eqref{unimockdef} conclude the proof of \Cref{sec:ummmflem1} (1).

	Next, using the definition of $g_2$ in \eqref{unimockdef} and the fact that \[\lim_{w\to0}\ \left(-\frac qw\right)_n w^n=q^\frac{n(n+1)}2,\] we have
	\begin{equation}\label{sec:ummmflem1eqn4}
	\left(1-\zeta\right)\left(1-\frac q\zeta\right)g_2(\zeta;q)=\lim_{w\to0}\sum_{n\ge 0}\frac{\left(-q,-\frac qw\right)_n}{\left(\zeta q,\frac{q^2}\zeta\right)_n}w^n.
	\end{equation}
		Now, letting $(a,b,c,d,e)=(q,\frac q\zeta,\zeta,-q,-\frac qw)$ in \Cref{sec:ummmflem1eqn5}, we get 
	\begin{align*}
	\sum_{n\ge 0}\frac{\left(-q,-\frac qw\right)_n}{\left(\zeta q,\frac{q^2}\zeta\right)_n}w^n=\frac{\left(-q,-wq,q\right)_{\infty}}{\left(\zeta q,\frac{q^2}\zeta,w\right)_{\infty}}\sum_{n\ge 0}\frac{\left(w,\frac q\zeta,\zeta\right)_n}{\left(q,-q,-wq\right)_n}q^n.
	\end{align*}
 Taking the limit $w\to0$ on both sides, and using \eqref{sec:ummmflem1eqn4}, we conclude the proof of \Cref{sec:ummmflem1} (2).\qedhere 
\end{proof}

We next prove \Cref{mmmfocase1lem2}.

\begin{proof}[Proof of \Cref{mmmfocase1lem2}]
	Using the definition of $\mathrm{op}^{[1]}_{r}(n)$, it is not hard to see that 
	\[
	\operatorname{OP}^{[1]}_{r}(q)=\sum_{n\ge 0} \frac{\left(-q^r;q^{2r}\right)^2_{n}}{\left(q^{2r};q^{2r}\right)_n}q^{2rn}.
	\]
Here, for $n\in \mathbb{N}$, the largest part is $2rn$, while $n=0$ accounts for the empty partition. Letting $q\mapsto q^{2r}$ and $\zeta=-q^r$ in \Cref{sec:ummmflem1} (1), we thus obtain the claim.\qedhere 
\end{proof} 

Now we establish \Cref{mmmfocase2lem1}.

\begin{proof}[Proof of \Cref{mmmfocase2lem1}]
	By definition of $p^{[2]}_{d,r}(n)$, we have
	\[
	P^{[2]}_{d,r}(q)=\sum_{n\ge 0}\frac{q^{dn}}{\left(q^d;q^d\right)_n\left(q^{dn+r},q^{dn+d-r};q^d\right)_{\infty}}.
	\]
Here, for $n\in \mathbb{N}$, the largest part divisible by $d$ is $dn$, while $n=0$ accounts for partitions with no parts divisible by $d$. Letting $q\mapsto q^{d}$ and $\zeta=q^r$ in \Cref{sec:ummmflem1} (1), we thus obtain the claim.\qedhere 
\end{proof}

Next, we prove \Cref{mmmfocase2lem2}.

\begin{proof}[Proof of \Cref{mmmfocase2lem2}]
	By definition of $p^{[2]}_{r}(n)$, we obtain
	\[P^{[2]}_{r}(q)=\sum_{n\ge 0}\frac{q^{2rn}}{\left(q^{2r};q^{2r}\right)_n\left(q^{2rn+r};q^{2r}\right)^2_{\infty}}.
	\]
	Letting $q\mapsto q^{2r}$ and $\zeta=q^r$ in \Cref{sec:ummmflem1} (1), we thus obtain the claim.\qedhere 
\end{proof}

 Next, we show the following elementary lemma which we require to prove Theorems \ref{mmmfecase1lem2}--\ref{mmmfecase2lem2}. 

\begin{lemma}\label{L:oddparts}
	For $k\in\mathbb{N}$ and $m\in\mathbb{N}_0$, $\frac{q^{km}}{(q^{2k};q^{2k})_m}$ is the generating function for partitions into exactly $m$ parts which are $\equiv k \pmod{2k}$.
\end{lemma}

\begin{proof}
Introduce a variable $w$ to record the number of parts. For each allowed part $k(2j+1)$ ($j\in \mathbb{N}_0$), we have
\[
\frac{1}{1-wq^{k\left(2j+1\right)}}=\sum_{r\ge 0} w^r q^{k(2j+1)r}.
\]
The term indexed by $r$ corresponds to taking $r$ copies of this part. Thus the coefficient of $w^m$ in
\[
\prod_{j\ge 0}\frac{1}{1-wq^{k\left(2j+1\right)}}=\frac{1}{\left(wq^k;q^{2k}\right)_\infty}
\]
is the generating function for partitions into exactly $m$ parts congruent to $k\!\pmod{2k}$.
 Letting $q\mapsto q^{2k}$ and $w\mapsto wq^k$ in \Cref{Euler}, we obtain
	\begin{equation*}
	\sum_{n\ge0}\frac{q^{kn}}{(q^{2k};q^{2k})_n}w^n=\frac{1}{\left(wq^k;q^{2k}\right)_\infty}.
	\end{equation*}
Comparing the coefficients of $w^m$ on both sides gives the claim.\qedhere 
\end{proof}

Now we show \Cref{mmmfecase1lem2}.

\begin{proof}[Proof of \Cref{mmmfecase1lem2}]
 Let $m\in\mathbb{N}_0$ denote the number of parts $\equiv 2r \pmod{4r}$ (including repetitions) of an overpartition counted by $\mathrm{op}^{[2]}_r(n)$. By (2), these parts are non-overlined. By \Cref{L:oddparts} with $k=2r$, they are generated by $\frac{q^{2rm}}{(q^{4r};q^{4r})_m}$. By (4), the parts $\equiv r\!\pmod{2r}$ lie in $\{r(2k-1):1\le k\le m\}$. Moreover, by (3), each of these parts appears at most once non-overlined and at most once overlined. Hence these parts are generated by
\[
\prod_{j=0}^{m-1}\left(1+q^{r\left(2j+1\right)}\right)^2=\left(-q^r;q^{2r}\right)_m^2.
\]
Multiplying these two generating functions and summing over $m$, we obtain
\[
\mathrm{OP}^{[2]}_r(q)=\sum_{m\ge0}\frac{\left(-q^r;q^{2r}\right)_m^2}{\left(q^{4r};q^{4r}\right)_m}q^{2rm}.
\]
Letting $q\mapsto q^{2r}$ and $\zeta=-q^r$ in \Cref{sec:ummmflem1} (2), we thus obtain the claim.\qedhere 
\end{proof}

Next, we prove \Cref{mmmfecase2lem1}.

\begin{proof}[Proof of \Cref{mmmfecase2lem1}]
Let $m\in\mathbb{N}_0$ be the number of parts $\equiv d \pmod{2d}$ (including repetitions) of a partition counted by $p^{[4]}_{d,r}(n)$. By \Cref{L:oddparts} with $k=d$, these parts are generated by $\frac{q^{dm}}{(q^{2d};q^{2d})_m}$. By (2), the parts $\equiv r\!\pmod d$ are at least $dm+r$ and the parts $\equiv -r \pmod d$ are at least $dm+d-r$. These parts are therefore generated by $\frac{1}{(q^{dm+r},q^{dm+d-r};q^d)_\infty}$. Multiplying these two generating functions and summing over $m$, we obtain
\[
P^{[4]}_{d,r}(q)=\sum_{m\ge0}\frac{q^{dm}}{\left(q^{2d};q^{2d}\right)_m\left(q^{dm+r},q^{dm+d-r};q^d\right)_\infty}.
\]
Letting $q\mapsto q^{d}$ and $\zeta=q^r$ in \Cref{sec:ummmflem1} (2), we thus obtain the claim.\qedhere 
\end{proof}

Finally, we establish \Cref{mmmfecase2lem2}.

\begin{proof}[Proof of \Cref{mmmfecase2lem2}]
 Let $m\in\mathbb{N}_0$ denote the number of parts $\equiv 2r \pmod{4r}$ (including repetitions) of a two-color partition counted by $p^{[3]}_r(n)$. \hspace{-0.05 cm}By (2), these parts appear in blue, so by \Cref{L:oddparts} with $k=2r$, they are generated by $\frac{q^{2rm}}{(q^{4r};q^{4r})_m}$. By (2) and (3), the parts $\equiv r\!\pmod{2r}$ are at least $2rm+r$. They are therefore generated by $\frac{1}{(q^{2rm+r};q^{2r})_\infty^2}$. Multiplying these two generating functions and summing over $m$, we obtain
\[
P^{[3]}_r(q)=\sum_{m\ge0}\frac{q^{2rm}}{\left(q^{4r};q^{4r}\right)_m\left(q^{2rm+r};q^{2r}\right)_\infty^2}.
\]	
Letting $q\mapsto q^{2r}$ and $\zeta=q^r$ in \Cref{sec:ummmflem1} (2), we thus obtain the claim.\qedhere 
\end{proof}

\section{Open problems}\label{prob}

We conclude this paper with the following problems for future research:
\begin{enumerate}
	\item Numerically,\footnote{Throughout, we verified all the congruences for $n\le 3000$.} we found the following congruences: 
	\begin{align*}
	{\rm op}^{[1]}_1(12n+8)&\equiv 0\pmod 2,\\
	{\rm op}^{[1]}_1(20n+j)&\equiv 0\pmod 2\ \text{with}\   j\in\{8,16\},\\
	{\rm op}^{[1]}_1(28n+j)&\equiv 0\pmod 2\ \text{with}\   j\in\{8,16,20\}.
	\end{align*}
	These lead to the following conjecture. Let $\left(\frac{\cdot}{\cdot}\right)$ denote the Legendre symbol.
	\begin{conjecture}
		For any odd prime $p$ and $k\in \mathbb{N}_0$, we have 
		\[
		2\mid {\rm op}^{[1]}_1\left(4pn+4k\right)\ \text{for all}\  n\in\mathbb{N}_0
		\Leftrightarrow \left(\frac{8k+1}{p}\right)=-1.
		\]
	\end{conjecture}
\end{enumerate}

\begin{enumerate}[resume]
	\item For ${\rm op}^{[2]}_1(n)$, numerically, we observe the following congruences:
	\begin{align*}
	{\rm op}^{[2]}_1(18n+j)&\equiv 0\pmod{4}\ \text{with}\  j\in\{8,14\},\\
	{\rm op}^{[2]}_1(50n+j)&\equiv 0\pmod{4}\ \text{with}\   j\in\{16,26,36,46\},\\
	{\rm op}^{[2]}_1(98n+j)&\equiv 0\pmod{4}\ \text{with}\   j\in\{26,40,54,68,82,96\}.
	\end{align*}
	These suggest the following general pattern.
	\begin{conjecture}
		Let $p$ be an odd prime with $p\not\equiv 1\!\pmod{8}$. Then for all $n\in\mathbb{N}_0$ and $1\le k\le p-1$, we have
		\[
		{\rm op}^{[2]}_1\!\left(2p^2n+\frac{p^2-1}{4}+2pk\right)\equiv 0\pmod{4}.
		\]
	\end{conjecture}
\end{enumerate}

\end{document}